\documentclass[12pt, a4paper]{article}
\usepackage[pagewise]{lineno}
\usepackage{graphicx}
\usepackage{amssymb}
\usepackage{latexsym, bm}
\usepackage{multicol}
\usepackage{indentfirst}
\usepackage{amssymb,amsfonts}
\usepackage{amsmath}
\usepackage{cases}
\usepackage[noadjust]{cite}
\usepackage[colorlinks=true,linkcolor=blue, citecolor=blue]{hyperref}
\newtheorem{theorem}{Theorem}
\newtheorem{lemma}{Lemma}

\newtheorem{definition}{Definition}

\newtheorem{remark}{Remark}

\usepackage{amssymb}
\numberwithin{equation}{section}
\numberwithin{theorem}{section}
\numberwithin{remark}{section}
\numberwithin{definition}{section}
\numberwithin{lemma}{section}
\numberwithin{corollary}{section}
\numberwithin{proposition}{section}

\title{Energy rigidity and weak--strong uniqueness for the
anisotropic Navier--Stokes equations in $\mathbb{R}^2_+$}
\author{Wenke Tan\footnote{tanwenkeybfq@163.com}\\
	{\small Key Laboratory of Computing and Stochastic Mathematics (Ministry of Education),}\\
	{\small School of Mathematics and Statistics, Hunan Normal University,}\\
	{\small Changsha, Hunan 410081, China}\\
 Fan Wu\footnote{wufan0319@yeah.net}\\
	{\small College of Science,	Jiangxi University of Water Resources and Electric Power,}\\
	{\small  Nanchang, Jiangxi 330099, China}\\
}

\date{}

\begin{document}
	\maketitle
	{\bf Abstract:} In this paper, we study the two-dimensional incompressible Navier--Stokes equations with only horizontal dissipation in the upper half-plane $\mathbb{R}^2_+$. The combination of the missing directional regularity and the influence of the boundary presents an essential difficulty for the pressure estimates.	
	We circumvent this difficulty by establishing an extension theorem for the weak formulation, based on the divergence-free condition and the homogeneous impermeability boundary condition. A key feature of this extension principle is that no boundary condition is required for the tangential velocity component. Consequently, every weak solution in the natural energy class on $\mathbb{R}^2_+$ admits a reflection extension that is a weak solution of the corresponding whole-space problem on $\mathbb{R}^2$. This extension allows us to transfer the energy rigidity and weak--strong uniqueness properties established by Demmel and Wiedemann \cite{E-W} to the half-plane setting. In particular, we prove the energy equality for weak solutions in the natural energy class and establish weak--strong uniqueness under the corresponding additional regularity assumptions.

	\medskip
	
	{\bf Mathematics Subject Classification (2020):} \  76W05, 76B03, 35Q35.
	\medskip
	
	{\bf Keywords:}  Anisotropic Navier--Stokes equations, Energy rigidity, Weak--strong uniqueness.
\section{Introduction}

In this paper, we consider the two-dimensional incompressible Navier--Stokes
equations with dissipation acting only in the horizontal direction
on the upper half-plane
$$
\mathbb R^2_+
:=\{x=(x_1,x_2)\in\mathbb R^2:\ x_2>0\}.
$$
More precisely, we study the initial-boundary value problem
\begin{equation}\label{ANS}
\left\{
\begin{aligned}
&\partial_t u-\partial_{x_1}^2u
+\nabla\cdot(u\otimes u)+\nabla P=0,\\
&\nabla\cdot u=0,
\end{aligned}
\right.
\end{equation}
in $(0,T)\times\mathbb R^2_+$, where
$u=(u_1,u_2)$ denotes the velocity field and $P$ is the pressure.
The system is supplemented with the impermeability boundary condition which is compatible with the divergence-free constraint
\begin{equation}\label{bc}
u_2(t,x_1,0)=0,
\end{equation}
and the initial condition
\begin{equation}\label{ic}
u(0,x)=u_0(x).
\end{equation}

The dissipation in \eqref{ANS} acts exclusively in the
$x_1$-direction. Consequently, the natural energy class provides
control of $u$ in $L^\infty(0,T;L^2(\mathbb R^2_+))$ and
of $\partial_{x_1}u$ in $L^2(0,T;L^2(\mathbb R^2_+))$,
but does not directly yield any control of the normal derivative
$\partial_{x_2}u_1$. This degeneracy creates a fundamental
difficulty in justifying the formal energy identity. The purpose of this paper is to
derive energy rigidity and weak--strong uniqueness for weak
solutions to \eqref{ANS}-\eqref{ic}. 
\subsection{Definitions and Notations}
Before presenting the main results, we give some notations and definitions.

\begin{definition}\label{def}
Let $T>0$. A function $u$ is called a weak solution to
\eqref{ANS}--\eqref{ic} on $[0,T)\times\mathbb{R}^2_+$
if the following conditions are satisfied:
\begin{enumerate}
    \item
    \[
    u\in L^\infty(0,T;L^2(\mathbb{R}^2_+)),
    \qquad
    \partial_{x_1}u\in L^2(0,T;L^2(\mathbb{R}^2_+)).
    \]

    \item
    For almost every $t\in[0,T)$,
    \[
    \int_{\mathbb{R}^2_+}u(t,x)\cdot\nabla\psi(x)\,dx=0
    \]
    for all $\psi\in C_c^\infty(\mathbb{R}^2_+)$.

    \item
    The following weak formulation holds:
    \begin{equation}\label{wf}
    \begin{aligned}
    &\int_0^T\int_{\mathbb{R}^2_+}
    \big(
    u\cdot\partial_t\phi
    -\partial_{x_1}u\cdot\partial_{x_1}\phi
    +u\otimes u:\nabla\phi
    \big)\,dx\,dt=
    \int_{\mathbb{R}^2_+}u_0(x)\cdot\phi(0,x)\,dx
    \end{aligned}
    \end{equation}
    for all $\phi\in C_c^\infty([0,T)\times\mathbb{R}^2_+;\mathbb{R}^2)$
    satisfying $\nabla\cdot\phi=0$.

    \item
    The boundary condition \eqref{bc} holds in the trace sense.
\end{enumerate}
\end{definition}

\begin{remark} To define weak solutions to system \eqref{ANS} on the whole space, it suffices to replace $\mathbb{R}^2_+$ by $\mathbb{R}^2$ in 1-3 of Definition \ref{def}. By the divergence-free condition, it follows 
\[
\partial_{x_2}u_2=-\partial_{x_1}u_1
\quad\text{in }\mathcal D'((0,T)\times\mathbb R^2_+).
\]
Since $\partial_{x_1}u_1\in
L^2(0,T;L^2(\mathbb R^2_+))$ and
$u_2\in L^\infty(0,T;L^2(\mathbb R^2_+))$,
we have
\[
u_2\in L^2(0,T;H^1(\mathbb R^2_+)).
\]
Moreover, its trace on $\partial\mathbb R^2_+$ vanishes
by the impermeability condition.
\end{remark}


Let $\Omega\subset \mathbb{R}^2$, we use $H^1(\Omega)$ to denote the space of measurable functions with the following norm
\begin{equation*}
\|f\|_{H^1(\Omega)}=(\int_{\Omega}(|f|^2+|\nabla f|^2)dx)^\frac{1}{2}.
\end{equation*}$L^q(0,T;L^p(\Omega))$ denotes the space of measurable functions with the following norm
\begin{equation*}
\|f\|_{L^q(0,T;L^p(\Omega)}
=
\begin{cases}
\displaystyle
\left(
\int_0^T
\left(
\int_{\Omega}|f(t,x)|^p\,dx
\right)^{\frac{q}{p}}
\,dt
\right)^{\frac{1}{q}},
& 1\leq q<\infty,\\[12pt]
\displaystyle
\operatorname*{ess\,sup}_{t\in(0,T)}
\|f(t,\cdot)\|_{L^p(\Omega)},
& q=\infty.
\end{cases}
\end{equation*}

Throughout this paper, $C$ is a positive constant and $C(A1,A2,...)$ is a positive
constant depending on $A1,A2,\cdots$.

\subsection{Main results}

We now state our main results.
\begin{theorem}[Energy rigidity]\label{main1}
Let $u$ be a weak solution to the problem \eqref{ANS}-\eqref{ic} in $[0,T)\times\mathbb{R}^2_+$. Then u has a unique representative in $C([0,T]:L^2(\mathbb{R}^2_+))$, which satisfies the anisotropic energy equality
\begin{align}\label{aee}
\int_{\mathbb{R}^2_+}|u(t,x)|^2dx+2\int_{0}^{t}\int_{\mathbb{R}^2_+}|\partial_{y_1} u(s,y)|^2dyds=\int_{\mathbb{R}^2_+}|u_0|^2dx
\end{align}for $0\leq t\leq T$.
\end{theorem}

Our second result, which follows from the anisotropic energy equality \eqref{aee}, shows the following weak--strong uniqueness principle.
\begin{theorem}[Weak--strong uniqueness]\label{main2}
 Let $u, v$ be weak solutions to \eqref{ANS}-\eqref{ic} with
$\partial_{x_2}v_1\in L^2(0,T;L^2(\mathbb{R}^2_+))$. If $u(0,x)=v(0,x)$ in $L^2(\mathbb{R}^2_+)$, then $u(t,x)=v(t,x)$ in $L^2(\mathbb{R}^2_+)$ for every $t\in [0,T)$.
\end{theorem}

\subsection{The known results and main ideas}\label{sec:known-results}
The energy equality for weak solutions is a fundamental
question in the mathematical theory of incompressible
fluid equations. For the classical Navier--Stokes equations,
Leray \cite{Leray} and Hopf \cite{Hopf} established the existence of global
weak solutions satisfying the energy inequality. Whether
this inequality can be upgraded to an equality is a
distinct issue at the weak-solution level, since the
formal cancellation of the nonlinear term requires
additional justification. This question has consequently
received considerable attention, leading to a variety
of sufficient conditions for energy equality.

For the 3d Navier--Stokes equations,
Lions \cite{Lions} and Prodi \cite{Prodi}  obtained the
classical energy equality criterion
$$
u\in L^4(0,T;L^4(\Omega)).
$$
Shinbrot \cite{Shinbrot} subsequently established
a broader class of sufficient conditions,
\[
u\in L^q(0,T;L^p(\Omega)),
\qquad
\frac{2}{q}+\frac{2}{p}=1,
\qquad p\geq4.
\]
Subsequent developments have produced further energy
equality criteria involving the velocity gradient,
the pressure, and various weak-in-time function spaces;
see, for example, \cite{Berselli,Cheskidov,Farwig,Kozono,Kukavica,Leslie1,Leslie,Luo,Shvydkoy,Yu}.
Related questions have also been investigated in the
context of Onsager's conjecture for the Euler equations.
For domains with boundaries, the analysis is further
complicated by the possible presence of boundary
contributions in the energy balance, making the
interplay between solution regularity and boundary
conditions particularly important \cite{Berselli,C-L,Farwig,Yu}.

A further challenge arises for fluid equations with
anisotropic dissipation. In contrast to the classical
Navier--Stokes equations, the viscous operator controls
derivatives only in selected spatial directions.
Consequently, the natural energy estimate does not
provide full spatial regularity, and the standard
justification of the energy identity cannot be applied
directly. The derivation of energy equality in this
setting therefore requires arguments that exploit
the interaction between the directional dissipation
and the underlying structure of the equations.

In this direction, Demmel and Wiedemann \cite{E-W} recently
studied the two-dimensional anisotropic Navier--Stokes
equations on the whole space $\mathbb R^2$:
$$
\left\{
\begin{aligned}
&\partial_t u+\nabla\cdot(u\otimes u)+\nabla p
=\partial_{x_1}^2u,\\
&\nabla\cdot u=0.
\end{aligned}
\right.
$$
They established an energy-rigidity result for weak
solutions in the natural energy class
$$
u\in L^\infty(0,T;L^2(\mathbb R^2)),
\qquad
\partial_{x_1}u\in L^2(0,T;L^2(\mathbb R^2)),
$$
without imposing any additional regularity in the
non-dissipative direction. Specifically, every weak
solution admits a unique representative in
$C([0,T];L^2(\mathbb R^2))$ that satisfies the
anisotropic energy equality for every $t\in[0,T]$.
In addition, they established a weak--strong uniqueness
criterion under the assumption that one of the two
solutions possesses the additional regularity
$$
\partial_{x_2}u_1\in L^2(0,T;L^2(\mathbb R^2)).
$$

Energy equality for anisotropic Navier--Stokes equations
has also been investigated in three dimensions.
In particular, Ye, Wei and Wang \cite{YWW} established
Lions--Shinbrot-type sufficient conditions for energy
equality in the presence of horizontal dissipation.
Their results address energy conservation under
additional integrability assumptions on weak solutions.
Together with the result of Demmel and Wiedemann \cite{E-W},
these works illustrate two complementary approaches
to energy equality for anisotropic fluid equations:
criteria based on additional integrability and
energy rigidity at the natural energy level.

The aim of the present paper is to establish the
corresponding energy-rigidity and weak-strong
uniqueness results for the two-dimensional
anisotropic Navier--Stokes equations on the upper
half-plane $\mathbb R^2_+$. The presence of the boundary 
leads to some essential difficulties. To illustrate the main obstructions and our ideas more clearly, we first review the essential ingredients of the analysis in \cite{E-W}. The first key ingredient in \cite{E-W} is the $L^2$ estimate for the pressure $p$. Taking the divergence of the equation \eqref{ANS}, we obtain the equation of pressure $p$
\begin{align*}
-\Delta p=\partial_{x_1}^2(u_1^2)+2\partial_{x_1}\partial_{x_2}(u_1u_2)+\partial_{x_2}^2(u_2^2),
\end{align*}
the key tool to obtain this estimate is a directional Riesz estimate for $u_1$,
\begin{align*}
\|R_1(u_1^2)\|_{L^2(\mathbb{R}^2)}\leq C\|u_1\|_{L^2(\mathbb{R}^2)}\|\partial_{x_1}u_1\|_{L^2(\mathbb{R}^2)}
\end{align*} where $\operatorname{R}_i$ is the Riesz transform in the i-th variable. However, this directional Riesz estimate, which relies on the Fourier transform, is not applicable in the half-space setting.
Even worse, because of the missing direction regularity of $u_1$, we have no suitable boundary conditions available for elliptic estimates of the pressure $p$. Another key ingredient in \cite{E-W} is the renormalization of $u_1$, which relies on the following commutator estimate.
\begin{align*}
\nabla\cdot (f^\epsilon v-(fv)^\epsilon)\to 0~strongly~in~L^1~for~\epsilon\to0
\end{align*}here $f$, $\partial_{x_1}f$, $v$ and $\partial_{x_1}v$ belong to $L^2(\mathbb{R}^2)$. Since this process is built on the mollification of $f$ and $fv$, while the half-space has no translation invariance in the vertical direction, there are also certain obstructions to the commutator estimate.



Our main contribution is to establish a
boundary-compatible extension principle for
weak solutions in the natural anisotropic energy
class, thereby circumventing the difficulties arising from boundary effects. The key step is to show that the weak
formulation remains valid for
test fields $\Phi\in C^\infty_c([0,T)\times\overline{\mathbb{R}^2_+})$ satisfying
\begin{align*}
\nabla\cdot\Psi=&0,\\
\Phi_2(t,x_1,0)=&0.
\end{align*}
We achieve this by employing
a shifted mollification in the normal direction,
which preserves the divergence-free constraint
while moving the support of the approximating
test fields away from the boundary. The resulting
normal derivative of the tangential component
contains a boundary-layer contribution that is
singular with respect to the mollification
parameter. We prove that this contribution
vanishes in the limit by using the zero trace
of $u_2$, its $L^2(0,T;H^1(\mathbb R^2_+))$
regularity, and Hardy's inequality. This argument
allows us to extend the admissible class of
test functions without imposing additional
regularity on the tangential velocity component.
We then define the extension $U=(U_1,U_2)$
to $\mathbb R^2$ by taking the even reflection
of the tangential component and the odd reflection
of the normal component:
$$
U_1(t,x_1,x_2)=u_1(t,x_1,|x_2|),
\qquad
U_2(t,x_1,x_2)
=\operatorname{sgn}(x_2)u_2(t,x_1,|x_2|).
$$
These parity choices are dictated by the
impermeability boundary condition and the
divergence-free constraint. In particular, the
zero trace of $u_2$ ensures that its odd extension
belongs to $H^1(\mathbb R^2)$ for almost every
time, while the even extension of $u_1$ preserves
the required $L^2$ and horizontal derivative
bounds. Combining the extended weak formulation
with the parity properties of the reflected
test fields, we establish that $U$ is a weak
solution of the corresponding anisotropic
Navier--Stokes equations on $\mathbb R^2$.

The two main theorems then follow from the
whole-space results of Demmel and Wiedemann \cite{E-W}.
For energy rigidity, the reflection doubles both
the kinetic energy and the horizontal dissipation,
so that the whole-space energy equality directly
yields the corresponding identity on the upper
half-plane. For weak--strong uniqueness, the
additional regularity
$$
\partial_{x_2}u_1\in
L^2(0,T;L^2(\mathbb R^2_+))
$$
is preserved by the even reflection of the
tangential component. The whole-space
weak--strong uniqueness theorem can therefore
be applied to the reflected solutions, and
the desired uniqueness result follows by
restriction to $\mathbb R^2_+$.

In summary, the present work transfers energy
rigidity and weak--strong uniqueness from the
whole plane to the upper half-plane through
a reflection principle compatible with the
natural energy regularity. The argument shows
that the impermeability condition, together
with incompressibility, provides sufficient
structure to handle the boundary without
requiring additional normal regularity of the
tangential velocity component. In particular,
energy rigidity holds throughout the natural
energy class, whereas weak--strong uniqueness
is obtained under the additional regularity
assumption stated above.


\section{Preliminary}
 
Our main tool in the proofs of Theorems \ref{main1} and
\ref{main2} is the following extension principle for weak solutions
to \eqref{ANS}--\eqref{ic}.

\begin{lemma}\label{etlm}
Let $u=(u_1,u_2)$ be a weak solution to
\eqref{ANS}--\eqref{ic} on
$[0,T)\times\mathbb{R}^2_+$.
Define $U=(U_1,U_2)$ on
$[0,T)\times\mathbb{R}^2$ by
\begin{equation*}
\begin{aligned}
U_1(t,x_1,x_2)
&=
\begin{cases}
u_1(t,x_1,x_2), & x_2\geq 0,\\
u_1(t,x_1,-x_2), & x_2<0,
\end{cases}\\[6pt]
U_2(t,x_1,x_2)
&=
\begin{cases}
u_2(t,x_1,x_2), & x_2\geq 0,\\
-u_2(t,x_1,-x_2), & x_2<0.
\end{cases}
\end{aligned}
\end{equation*}
Then $U$ is a weak solution to \eqref{ANS} on
$[0,T]\times\mathbb{R}^2$.
\end{lemma}
\begin{remark} Once the whole-space weak solution $U$ has been constructed, we can define the corresponding pressure $P$ by
\begin{align*}
P=\operatorname{R}_1^2(U_1^2)+2\operatorname{R}_1\operatorname{R}_2(U_1U_2)+\operatorname{R}_2^2(U_2^2). 
\end{align*}Then, the $L^2$ estimate for $P$ can be established by the same arguments developed in \cite{E-W}. Since $U|_{\mathbb{R}^2_+}=u$, we have
\begin{align*}
\Delta p=\Delta P~on~\mathbb{R}^2_+.
\end{align*}Thus, we establish an $L^2$ estimate for the original pressure $p$ modulo a harmonic function. We also wish to point out that our extension theorem remains valid in higher dimensions and is of independent interest.\end{remark}

{\bf Proof.}
\textbf{Step 1:  Extension of the test function.} We will show that the half-plane weak formulation \eqref{wf} remains valid for $\Phi(t,x_1,x_2)\in C^\infty_c([0,T)\times\overline{\mathbb{R}^2_+})$ satisfying $\nabla\cdot\Phi=0$ and $\Phi_2(t,x_1,0)=0$. No boundary condition is imposed on $\Phi_1$.

Let $\rho\in C_c^\infty(-1,1)$ be an even function satisfying
$$
\int_{\mathbb{R}}\rho(x)\,dx=1.
$$
We define
\begin{equation}\label{mt}
\Phi^\epsilon(t,x_1,x_2)
=
\int_0^\infty
\frac{1}{\epsilon}
\rho\left(\frac{x_2-y-2\epsilon}{\epsilon}\right)
\Phi(t,x_1,y)\,dy.
\end{equation}
By the definition of $\Phi^\epsilon$ and the properties of $\rho$,
we have
\begin{align}
\operatorname{supp}\Phi^\epsilon
&\subset
\left\{
(t,x_1,x_2)\in[0,T)\times\mathbb{R}^2_+
\,:\,x_2\geq\epsilon
\right\},
\label{sup}\\
\nabla\cdot\Phi^\epsilon&=0.
\label{df1}
\end{align}
Consequently, taking $\Phi^\epsilon$ as a test function in \eqref{wf}, we obtain
\begin{align}\label{apwf}
\int_{0}^{T}\int_{\mathbb{R}^2_+}u\cdot \partial_t\Phi^\epsilon-\partial_{x_1} u\cdot\partial_{x_1}\Phi^\epsilon+u\otimes u:\nabla\Phi^\epsilon dxdt=\int_{\mathbb{R}^2_+}u_0(x)\Phi^\epsilon(0,x)dx.\end{align}
We next establish several properties of $\Phi^\epsilon$
that will be used to pass to the limit in \eqref{apwf}
as $\epsilon\to 0$. By the smoothness of $\Phi$ and the
properties of the mollifier $\rho$, we have
\begin{align}
\partial_t\Phi^\epsilon
&\longrightarrow \partial_t\Phi
\quad\text{in }
C([0,T)\times\overline{\mathbb{R}^2_+}),
\label{a1}\\
\partial_{x_1}\Phi^\epsilon
&\longrightarrow \partial_{x_1}\Phi
\quad\text{in }
C([0,T)\times\overline{\mathbb{R}^2_+})
\label{a2}
\end{align}
as $\epsilon\to 0$.
For the normal derivative, an integration by parts yields
\begin{align*}
\partial_{x_2}\Phi^\epsilon_j(t,x_1,x_2)=&\int_{0}^{\infty}\frac{1}{\epsilon}\partial_{x_2}\rho(\frac{x_2-y-2\epsilon}{\epsilon})\Phi(t,x_1,y)dy\\
=&-\int_{0}^{\infty}\frac{1}{\epsilon}\partial_{y}\rho(\frac{x_2-y-2\epsilon}{\epsilon})\Phi_j(t,x_1,y)dy\\
=&\frac{1}{\epsilon}\rho(\frac{x_2-2\epsilon}{\epsilon})\Phi_j(t,x_1,0)+\int_{0}^{\infty}\frac{1}{\epsilon}\rho(\frac{x_2-y-2\epsilon}{\epsilon})\partial_y\Phi_j(t,x_1,y)dy
\end{align*}where $j\in\{1,2\}$.
Using the boundary condition $\Phi_2(t,x_1,0)=0$, we see that
\begin{align*}
\partial_{x_2}\Phi^\epsilon_2=\int_{0}^{\infty}\frac{1}{\epsilon}\rho(\frac{x_2-y-2\epsilon}{\epsilon})\partial_y\Phi_j(t,x_1,y)dy,\end{align*} and thus
\begin{align}
\lim_{\epsilon\to 0}\partial_{x_2}\Phi^\epsilon_2=\partial_{x_2}\Phi_2~in~C([0,T)\times\overline{\mathbb{R}^2_+}).\label{a3}
\end{align} Moreover, we have
\begin{align}
\partial_{x_2}\Phi^\epsilon_1
&=
\frac{1}{\epsilon}
\rho\left(\frac{x_2-2\epsilon}{\epsilon}\right)
\Phi_1(t,x_1,0)+
\int_0^\infty
\frac{1}{\epsilon}
\rho\left(\frac{x_2-y-2\epsilon}{\epsilon}\right)
\partial_y\Phi_1(t,x_1,y)\,dy
\label{a4}\\
&=:K_1(t,x_1,x_2)+K_2(t,x_1,x_2).
\nonumber
\end{align}
Here, $K_1$ denotes the singular boundary term, given by
\[
K_1(t,x_1,x_2)
=
\frac{1}{\epsilon}
\rho\left(\frac{x_2-2\epsilon}{\epsilon}\right)
\Phi_1(t,x_1,0),
\]
while $K_2$ is given by
\[
K_2(t,x_1,x_2)
=
\int_0^\infty
\frac{1}{\epsilon}
\rho\left(\frac{x_2-y-2\epsilon}{\epsilon}\right)
\partial_y\Phi_1(t,x_1,y)\,dy.
\]
By the support, boundedness, and normalization properties of $\rho$, we have
\begin{align}
&\operatorname{supp}K_1
\subset
\left\{
(t,x)\in[0,T)\times\overline{\mathbb{R}^2_+}
\,:\,
\epsilon\leq x_2\leq 3\epsilon
\right\},
\label{a5}\\
&|K_1(t,x_1,x_2)|
\leq \frac{C}{\epsilon},
\label{a6}\\
&K_2
\longrightarrow \partial_{x_2}\Phi_1
\quad\text{as }\epsilon\to 0.
\label{a7}
\end{align}
We now show the convergence of \eqref{apwf} term by term as $\epsilon\to 0$. Since $\Phi\in C^\infty_c([0,T)\times\overline{\mathbb{R}^2_+})$, we immediately deduce that
\begin{align}\label{c0}
\lim_{\epsilon\to 0}\int_{\mathbb{R}^2_+}u_0(x)\Phi^\epsilon(0,x)dx=\int_{\mathbb{R}^2_+}u_0(x)\Phi(0,x)dx.\end{align} For the time-derivative term, by using \eqref{a1}, it follows that
\begin{align}\label{c1}
\lim_{\epsilon\to0}\int_{0}^{T}\int_{\mathbb{R}^2_+}u\cdot\partial_{t}\Phi^\epsilon dxdt=\int_{\mathbb{R}^2_+}u\cdot\partial_{t}\Phi dxdt.\end{align} For the anisotropic dissipation term, we use \eqref{a2} to deduce that
\begin{align}\label{c2}
\lim_{\epsilon\to0}\int_{0}^{T}\int_{\mathbb{R}^2_+}\partial_{x_1}u\cdot\partial_{x_1}\Phi^\epsilon dxdt=\int_{\mathbb{R}^2_+}\partial_{x_1}u\cdot\partial_{x_1}\Phi dxdt.\end{align}

Next, we consider the convergence of the nonlinear term. By direct computations, we see that
\begin{align*}
u\otimes u:\nabla\Phi^\epsilon=u_1^2\partial_{x_1}\Phi^\epsilon_1+u_1u_2\partial_{x_2}\Phi^\epsilon_1+u_2u_1\partial_{x_1}\Phi^\epsilon_2+u_2^2\partial_{x_2}\Phi^\epsilon_2.
\end{align*}Since $u\in L^\infty(0,T;L^2(\mathbb{R}^2_+))$, from \eqref{a2} and \eqref{a3}, we obtain that
\begin{align}
&\lim_{\epsilon\to 0}\int_{0}^{T}\int_{\mathbb{R}^2_+}(u_1^2\partial_{x_1}\Phi^\epsilon_1+u_2u_1\partial_{x_1}\Phi^\epsilon_2+u_2^2\partial_{x_2}\Phi^\epsilon_2)dxdt\label{c3}\\
=&\int_{0}^{T}\int_{\mathbb{R}^2_+}(u_1^2\partial_{x_1}\Phi_1+u_2u_1\partial_{x_1}\Phi_2+u_2^2\partial_{x_2}\Phi_2)dxdt.\notag\end{align}
Then, using \eqref{a4}, we get
\begin{align}
\int_{0}^{T}\int_{\mathbb{R}^2_+}u_1u_2\partial_{x_2}\Phi^\epsilon_1dxdt=\int_{0}^{T}\int_{\mathbb{R}^2_+}u_1u_2K_1dxdt+\int_{0}^{T}\int_{\mathbb{R}^2_+}u_1u_2K_2dxdt\label{c4}.\end{align}From \eqref{a7}, we have
\begin{align}
\lim_{\epsilon\to0}\int_{0}^{T}\int_{\mathbb{R}^2_+}u_1u_2K_2dxdt=\int_{0}^{T}\int_{\mathbb{R}^2_+}u_1u_2\partial_{x_2}\Phi_1dxdt.\label{c5}\end{align}
It remains to estimate the contribution from $K_1$.
By \eqref{a5} and \eqref{a6}, we obtain
\begin{align*}
&\left|
\int_0^T\int_{\mathbb{R}^2_+}
u_1u_2K_1\,dx\,dt
\right|\\
\leq&
C\int_0^T\int_{\mathbb{R}}
\int_\epsilon^{3\epsilon}
|u_1|\frac{|u_2|}{\epsilon}
\,dx_2\,dx_1\,dt\\
\leq&
C\int_0^T\int_{\mathbb{R}}
\int_0^{3\epsilon}
|u_1|
\left|\frac{u_2}{x_2}\right|
\,dx_2\,dx_1\,dt.
\end{align*}

Since
$u_1\in L^\infty(0,T;L^2(\mathbb{R}^2_+))$
and
$u_2\in L^2(0,T;H^1(\mathbb{R}^2_+))$
with zero trace on $x_2=0$, the Hardy inequality and
the Cauchy--Schwarz inequality yield
\begin{align*}
&\int_0^T\int_{\mathbb{R}}\int_0^{3\epsilon}
|u_1|\left|\frac{u_2}{x_2}\right|
\,dx_2\,dx_1\,dt\\
\leq&
\left(
\int_0^T\int_{\mathbb{R}}\int_0^{3\epsilon}
|u_1|^2\,dx_2\,dx_1\,dt
\right)^{\frac12}
\left(
\int_0^T\int_{\mathbb{R}}\int_0^{3\epsilon}
|\partial_{x_2}u_2|^2\,dx_2\,dx_1\,dt
\right)^{\frac12},
\end{align*}
we thus get, from the absolute continuity of their integrals, that
\begin{align}\label{c6}
\lim_{\epsilon \to 0}\int_{0}^{T}\int_{\mathbb{R}^2_+}u_1u_2K_1dxdt=0.\end{align}
Substituting \eqref{c5} and \eqref{c6} into \eqref{c4} yields
\begin{align}\label{c7}
\lim_{\epsilon\to 0}\int_{0}^{T}\int_{\mathbb{R}^2_+}u_1u_2\partial_{x_2}\Phi^\epsilon_1dxdt
=\int_{0}^{T}\int_{\mathbb{R}^2_+}u_1u_2\partial_{x_2}\Phi_1dxdt.\end{align}
Combining \eqref{c0}-\eqref{c3} with \eqref{c7}, we deduce from \eqref{apwf} that
\begin{align}\label{ewf}
\int_{0}^{T}\int_{\mathbb{R}^2_+}u\cdot \partial_t\Phi-\partial_{x_1} u\cdot\partial_{x_1}\Phi+u\otimes u:\nabla\Phi dxdt=\int_{\mathbb{R}^2_+}u_0(x)\Phi(0,x)dx\end{align} for $\Phi(t,x_1,x_2)\in C^\infty_c([0,T)\times\overline{\mathbb{R}^2_+})$ satisfying $\nabla\cdot\Phi=0$ and $\Phi_2(t,x_1,0)=0$.

\textbf{Step 2: The extended weak solution on $\mathbb{R}^2$.}
We define $U=(U_1,U_2)$ on $\mathbb{R}^2$ by
\begin{equation}\label{eu1}
U_1(t,x_1,x_2)
=
\begin{cases}
u_1(t,x_1,x_2), & x_2\geq 0,\\
u_1(t,x_1,-x_2), & x_2<0,
\end{cases}
\end{equation}
and
\begin{equation}\label{eu2}
U_2(t,x_1,x_2)
=
\begin{cases}
u_2(t,x_1,x_2), & x_2\geq 0,\\
-u_2(t,x_1,-x_2), & x_2<0.
\end{cases}
\end{equation}
By construction, we have
\[
U\in L^\infty(0,T;L^2(\mathbb{R}^2)),
\qquad
\partial_{x_1}U\in L^2(0,T;L^2(\mathbb{R}^2)).
\]
Moreover, the following norm identities hold for almost every $t\in(0,T)$:
\begin{align}
\|U(t)\|_{L^2(\mathbb{R}^2)}^2
&=2\|u(t)\|_{L^2(\mathbb{R}^2_+)}^2,
\label{e1}\\
\|\partial_{x_1}U(t)\|_{L^2(\mathbb{R}^2)}^2
&=2\|\partial_{x_1}u(t)\|_{L^2(\mathbb{R}^2_+)}^2.
\label{e2}
\end{align}
Since $u_2(t)\in H^1(\mathbb{R}^2_+)$ and has zero trace on $x_2=0$, $U_2$ is the odd extension of $u_2$ with respect to $x_2$, it follows that $U_2(t)\in H^1(\mathbb{R}^2)$. Moreover,  $U_1$ is the even extension of $u_1$, and its tangential derivative $\partial_{x_1}U_1$ is the even extension of $\partial_{x_1}u_1$. A direct computation implies
\begin{align}\label{df}
\partial_{x_1}U_1+\partial_{x_2} U_2=0.
\end{align}
We will prove that $U(t,x)=(U_1(t,x),U_2(t,x))$ is a weak solution to \eqref{ANS} on $\mathbb{R}^2$. Let $\psi(t,x)\in C^\infty_c([0,T)\times\mathbb{R}^2)$ obeying $\nabla\cdot\psi=0$. By the definition of $U_1$, $U_2$ and the change of variables,
we have that
\begin{align}\label{ewf1}
&\int_{0}^{T}\int_{\mathbb{R}^2}U_1\partial_t\psi_1-\partial_{x_1}U_1\partial_{x_1}\psi_1+U_1^2\partial_{x_1}\psi_1+U_1U_2\partial_{x_2}\psi_1dxdt\\
=&\int_{0}^{T}\int_{\mathbb{R}^2_+}u_1\partial_t\psi_1-\partial_{x_1}u_1\partial_{x_1}\psi_1+u_1^2\partial_{x_1}\psi_1+u_1u_2\partial_{x_2}\psi_1dxdt\notag\\
+&\int_{0}^{T}\int_{\mathbb{R}^2_{-}}u_1(\cdot,-x_2)\partial_t\psi_1-\partial_{x_1}u_1(\cdot,-x_2)\partial_{x_1}\psi_1+u_1^2(\cdot,-x_2)\partial_{x_1}\psi_1-(u_1u_2)(\cdot,-x_2)\partial_{x_2}\psi_1dxdt\notag\\
=&\int_{0}^{T}\int_{\mathbb{R}^2_+}u_1\partial_t\Psi_1-\partial_{x_1}u_1\partial_{x_1}\Psi_1+u_1^2\partial_{x_1}\Psi_1+u_1u_2\partial_{x_2}\Psi_1dxdt,\notag\end{align} here the function $\Psi_1$ is defined as
\begin{align}\label{Psi1}
\Psi_1(t,x)=\psi_1(t,x_1,x_2)+\psi_1(t,x_1,-x_2).
\end{align} Similarly,
\begin{align}\label{ewf2}
&\int_{0}^{T}\int_{\mathbb{R}^2}U_2\partial_t\psi_2-\partial_{x_1}U_2\partial_{x_1}\psi_2+U_2U_1\partial_{x_1}\psi_2+U_2^2\partial_{x_2}\psi_2dxdt\\
=&\int_{0}^{T}\int_{\mathbb{R}^2_+}u_2\partial_t\psi_2-\partial_{x_1}u_2\partial_{x_1}\psi_2+u_2u_1\partial_{x_1}\psi_2+u_2^2\partial_{x_2}\psi_2dxdt\notag\\
+&\int_{0}^{T}\int_{\mathbb{R}^2_{-}}-u_2(\cdot,-x_2)\partial_t\psi_2+\partial_{x_1}u_2(\cdot,-x_2)\partial_{x_1}\psi_2-(u_2u_1)(\cdot,-x_2)\partial_{x_1}\psi_2+u_2^2(\cdot,-x_2)\partial_{x_2}\psi_2dxdt\notag\\
=&\int_{0}^{T}\int_{\mathbb{R}^2_+}u_2\partial_t\Psi_2-\partial_{x_1}u_2\partial_{x_1}\Psi_2+u_2u_1\partial_{x_1}\Psi_2+u_2^2\partial_{x_2}\Psi_2dxdt,\notag\end{align}here $\Psi_2$ is defined as
\begin{align}\label{Psi2}
\Psi_2(t,x)=\psi_2(t,x_1,x_2)-\psi_2(t,x_1,-x_2).
\end{align}
Let $\Psi(t,x)=(\Psi_1(t,x),\Psi_2(t,x))$. It holds that
\begin{align}
&\Psi\in C^\infty_c([0,T)\times\overline{\mathbb{R}^2_+}),\\
&\Psi_2(t,x_1,0)=0,\\
\nabla \cdot\Psi=&\partial_{x_1}\psi_1(t,x_1,x_2)+\partial_{x_1}\psi_1(t,x_1,-x_2)\\
&+\partial_{x_2}\psi_2(t,x_1,x_2)-\partial_{x_2}\psi_2(t,x_1,-x_2)\notag\\
=&(\nabla\cdot\psi)(t,x_1,x_2)+(\nabla\cdot\psi)(t,x_1,-x_2)\notag\\
=&0.\notag\end{align}
By using the conclusion established in {\bf Step 1}, we deduce from \eqref{ewf}, \eqref{ewf1} and \eqref{ewf2} that
\begin{align}\label{lwf1}
&\int_{0}^{T}\int_{\mathbb{R}^2}U_1\partial_t\psi_1-\partial_{x_1}U_1\partial_{x_1}\psi_1+U_1^2\partial_{x_1}\psi_1+U_1U_2\partial_{x_2}\psi_1dxdt\\
=&\int_{0}^{T}\int_{\mathbb{R}^2_+}u_1\partial_t\Psi_1-\partial_{x_1}u_1\partial_{x_1}\Psi_1+u_1^2\partial_{x_1}\Psi_1+u_1u_2\partial_{x_2}\Psi_1dxdt,\notag\\
=&\int_{\mathbb{R}^2_+}u_1(0,x)(\psi_1(0,x_1,x_2)+\psi_1(0,x_1,-x_2))dx\notag\\
=&\int_{\mathbb{R}^2}U_1(0,x)\psi_1(0,x)dx\notag,
\end{align}
\begin{align}\label{lwf2}
&\int_{0}^{T}\int_{\mathbb{R}^2}U_2\partial_t\psi_2-\partial_{x_1}U_2\partial_{x_1}\psi_2+U_2U_1\partial_{x_1}\psi_2+U_2^2\partial_{x_2}\psi_2dxdt\\
=&\int_{0}^{T}\int_{\mathbb{R}^2_+}u_2\partial_t\Psi_2-\partial_{x_1}u_2\partial_{x_1}\Psi_2+u_2u_1\partial_{x_1}\Psi_2+u_2^2\partial_{x_2}\Psi_2dxdt,\notag\\
=&\int_{\mathbb{R}^2_+}u_2(0,x)(\psi_2(0,x_1,x_2)-\psi_2(t,x_1,x_2))dx\notag\\
=&\int_{\mathbb{R}^2}U_2(0,x)\psi_2(0,x)dx.\notag
\end{align}Combining \eqref{e1}, \eqref{e2}, \eqref{df}, \eqref{lwf1}, \eqref{lwf2}, we conclude that
the function $U=(U_1,U_2)$ is a weak solution to \eqref{ANS} with initial data $U(0,x)$.\qquad $\hfill\Box$

\section{Proofs of main theorems}
\par

\textbf{The proof of Theorem \ref{main1}.} Applying the energy rigidity theorem (Theorem 1.2 in \cite{E-W}) for weak solutions on $\mathbb{R}^2$, we obtain that
\begin{align}
\int_{\mathbb{R}^2}|U(t,x)|^2dx+2\int_{0}^{t}\int_{\mathbb{R}^2}|\partial_{x_1}U(s,x)|^2dxds=\int_{\mathbb{R}^2}|U(0,x)|^2dx.
\end{align}Substituting \eqref{e1} and \eqref{e2} into the above identity yields
\begin{align}
\int_{\mathbb{R}^2_+}|u(t,x)|^2dx+2\int_{0}^{t}\int_{\mathbb{R}^2_+}|\partial_{x_1}u(s,x)|^2dxds=\int_{\mathbb{R}^2_+}|u(0,x)|^2dx,\end{align} we obtain the energy rigidity for weak solutions on $\mathbb{R}^2_+$.\qquad $\hfill\Box$

\textbf{The proof of Theorem \ref{main2}.} Assume $U$ and  $V$ are the extended weak solutions of $u$ and $v$ following Lemma \ref{etlm}. We have seen that  $U$ and $V$ are weak solutions on $\mathbb{R}^2$ with the same initial data. Since $\partial_{x_2}v_1\in L^2([0,T);\mathbb{R}^2_+)$ and $V_1$ is the even extension of $v_1$ with respect to $x_2$, it follows that $\partial_{x_2}V_1\in L^2([0,T;L^2(\mathbb{R}^2)))$. Applying the weak-strong uniqueness theorem of Demmel and Wiedemann (Theorem 1.3 in \cite{E-W}), we get that $U(t)=V(t)$ in $L^2(\mathbb{R}^2)$. From the definitions of $U$ and $V$, it follows $u(t)=v(t)$ in $L^2(\mathbb{R}^2_+)$. \qquad $\hfill\Box$

\subsection*{Acknowledgments} W. Tan was supported by the Construct Program of the Key Discipline in Hunan Province and the National Natural Science Foundation of China (11871209);
F. Wu was supported by Jiangxi Provincial Natural
	Science Foundation (20252BAC240160).
	
\section*{Conflict of Interest}
The authors have no conflicts to disclose.

\section*{Data Availability}
Data sharing is not applicable to this article as no new data were created or analyzed in this study.

\end{document}